\documentclass[12pt]{article}
\usepackage{amsfonts,amsthm,amsmath,amssymb}
\usepackage{latexsym,graphicx,fullpage}
\usepackage{hypernat,titlesec,titletoc}

\title{{The Minimum Stretch Spanning Tree Problem\\
for Typical Graphs}
\thanks{Supported by NSFC (61373106) and 973 Program of China
(2010CB328101).}}
\author{{Lan Lin$^{1}$, \ \ \ Yixun Lin$^{2}$}\thanks{Corresponding Author. E-mail address:
linlan@tongji.edu.cn, linyixun@zzu.edu.cn}\\
{\small $^{1}$ School of Electronics and Information Engineering,
Tongji University,}\\ {\small Shanghai 200092, China. }\\
{\small $^{2}$ School of Mathematics and Statistics, Zhengzhou
University, Zhengzhou 450001, China.}}

\date{}
\begin{document}
\maketitle
\renewcommand\baselinestretch{1.0}

{\noindent\bf Abstract} \ With applications in communication
networks, the minimum stretch spanning tree problem is to find a
spanning tree $T$ of a graph $G$ such that the maximum distance in
$T$ between two adjacent vertices is minimized. The problem has been
proved to be NP-hard and fixed-parameter polynomial algorithms have
been obtained for some special families of graphs. In this paper, we
concentrate on the optimality characterizations for typical classes
of graphs. We determine the exact results for the Petersen graph,
the complete $k$-partite graphs, split graphs, generalized convex
graphs, and several planar grids, including rectangular grids,
triangular grids, and triangulated-rectangular grids.


{\noindent\bf 2000 MR Subject Classification} \  90C27, 05C05

\section{Introduction }
\hspace*{0.5cm} Since Peleg et al.\cite{Peleg89} in 1989, a series
of tree spanner problems arise in connection with applications in
distribution systems and communication networks (see survey
\cite{Lieb08}). A basic decision version of the tree spanner
problems for a graph $G$ is as follows: For a given integer $k$, is
there a spanning tree $T$ of $G$ (called a tree $k$-spanner) such
that the distance in $T$ between every pair of vertices is at most
$k$ times their distance in $G$? The corresponding optimization
version of the problem is to find the minimum $k$ such that there
exists a tree $k$-spanner of $G$. This spanning tree optimization
problem is referred to as {\it the minimum stretch spanning tree
problem} and MSST for short \cite{Brand04,Cai95,Fekete01,Madan96}.

We formulate the problem formally. Let $G$ be a simple connected
graph with vertex set $V(G)$ and edge set $E(G)$. Given a spanning
tree $T$ of $G$, for $uv\in E(G)$, let $d_T(u,v)$ denote the
distance between $u$ and $v$ in $T$, that is the length of the
unique $u$-$v$-path in $T$. Then the {\it stretch} of a spanning
tree $T$ is defined by
\begin{equation}
\sigma_T(G,T):=\max_{uv\in E(G)}\,d_T(u,v).
\end{equation}
Furthermore, the minimum stretch spanning tree problem is to
determine
\begin{equation}
\sigma_T(G):=\min \{\sigma_T(G,T): T \mbox{ is a spanning tree of }
G\}.
\end{equation}
This gives rise to a graph invariant $\sigma_T(G)$, called {\it the
tree-stretch} of $G$. Here, we follow the notation $\sigma_T(G)$ in
\cite{Fekete01}.

For an edge $e=uv$ not in $T$, the unique cycle in $T+e$ is called
the {\it fundamental cycle} with respect to $e$. So, the above
problem is equivalent to finding a spanning tree such that the
length of a maximum fundamental cycle is minimized, where the
tree-stretch $\sigma_T(G)$ is one less than the length of this
cycle. This is precisely {\it the shortest maximal fundamental cycle
problem} proposed by Galbiati \cite{Galb03}. As is well known, all
fundamental cycles with respect to a spanning tree $T$ constitute a
basis of the {\it cycle space} of $G$ \cite{Bondy08}. Thus we have
an optimal basis problem in the cycle space.

In the dual point of view, for each $e\in T$, the edge-cut between
two components of $T-e$ is a fundamental edge-cut (cocycle). Let
$X_e$ be the vertex set of one of these components. Write $\partial
(X_e):=\{uv\in E(G):u\in X_e, v\notin X_e\}$. Then $\partial (X_e)$
is the {\it fundamental edge-cut} with respect to $e$, and
$|\partial (X_e)|$ is called the {\it congestion} of edge $e$. The
minimum congestion spanning tree problem, proposed by Ostrovskii
\cite{Ostrov04} in 2004, is to determine
$$c_T(G):=\min \{\max_{e\in T}\,|\partial(X_e)|:T \mbox{ is a spanning tree of } G\}.$$
This graph invariant $c_T(G)$ is called {\it the tree-congestion} of
$G$.

Admittedly, the tree-congestion $c_T(G)$ is a variant of the
cutwidth $c(G)$ of $G$ and the tree-stretch $\sigma_T(G)$ is a
variant of the bandwidth $B(G)$ of $G$  (see surveys
\cite{Chung88,Diaz02}). In the circuit layout of VLSI designs and
network communication, the quality of an embedding is usually
evaluated by two parameters, namely, the dilation and the
congestion. The dilation motivates the bandwidth problem and the
congestion leads to the cutwidth problem.

So far the main concern of the tree spanner problems is in the
algorithmic aspects, including the NP-hardness \cite{Brand04,
Brand07,Cai95,Fekete01,Galb03}, the fixed-parameter polynomial
algorithms \cite{Brand04,Brand07,Fekete01,Fomin11}, and the
approximability \cite{Galb03}. Moreover, for the characterization
problem, it is known that determining $\sigma_T(G)\leq 2$ is
polynomially solvable \cite{Cai95}, while determining $\sigma_T\leq
k$ for $k\geq 4$ is NP-complete. A long-standing open problem is to
characterize $\sigma_T(G)=3$. In this respect, it is significant to
determine exact value of $\sigma_T(G)$ for typical classes of
graphs.

The minimum congestion spanning tree problem has been studied
extensively in the literature. On the complexity aspect, the
NP-hardness even for chain graphs or split graphs was shown in
\cite{Okamoto11}. Linear time algorithms for fixed parameter $k$ and
for planar graphs, bounded-degree graphs and treewidth bounded
graphs were presented in \cite{Bodlae12}. Additionally, determining
the exact values of $c_T(G)$ for special graphs has found an
increasing interest during the last decade, for example:\\
{\indent} $\bullet$ The complete graphs $K_n$, the complete
bipartite graphs $K_{m,n}$, and the planar grids $P_m\times P_n$
\cite{Caste09,Hruska08}.\\
{\indent} $\bullet$ The complete $k$-partite graphs
$K_{n_1,n_2,\ldots,n_k}$ and the torus grids $C_m\times C_n$
\cite{Caste09,Kozawa09}. \\
{\indent} $\bullet$ The triangular grids $T_n$ \cite{Ostrov10}. \\
{\indent} $\bullet$ The $k$-outerplanar graphs \cite{Bodlae11}.

Motivated by the above results on $c_T(G)$, our goal is to
investigate the dual invariant $\sigma_T(G)$ for some basic families
of graphs. The main results are parallel to those for $c_T(G)$.

The remaining of the paper is organized as follows. In Section 2, we
present a basic lower bound by using the girth and derive the exact
results for $K_n$, $C_n$, $K_{m,n}$, the Petersen graph, etc. In
Section 3, we characterize $K_{n_1,n_2,\ldots, n_k}$, split graphs
and generalized convex graphs. Section 4 is devoted to the exact
representations for a class of plane graphs, including rectangular
grids $P_m\times P_n$, triangular grids $T_n$, and
triangulated-rectangular grids $T_{m,n}$.

\section{Elementary properties}
\hspace*{0.5cm} We shall follow the graph-theoretic terminology and
notation of \cite{Bondy08}. Let $G$ be a simple connected graph on
$n$ vertices with vertex set $V(G)$ and edge set $E(G)$. For a
subset $S\subseteq V(G)$, the {\it neighbor set} of $S$ is defined
by $N_G(S):=\{v\in V(G)\setminus S: u\in S,uv\in E(G)\}$. We
abbreviate $N_G(\{v\})$ to $N_G(v)$ for a vertex $v\in V(G)$. For
$S\subseteq V(G)$, we denote by $G[S]$ the subgraph induced by $S$.
For an edge $e\in E(G)$, denote by $G-e$ the graph obtained from $G$
by deletion of $e$. For an edge $e$ not in $E(G)$, denote by $G+e$
the graph obtained from $G$ by addition of $e$.

Let $T$ be a spanning tree of $G$. As usual, the spanning tree $T$
is regarded as a set of edges. The {\it cotree} $\bar T$ of $T$ is
defined as the complement of $T$ in $E(G)$, namely  $\bar T=
E(G)\setminus T$. For each $e\in \bar T$, the unique cycle in $T+e$
is a fundamental cycle, determined by the cotree edge $e$. The
tree-stretch $\sigma_T(G)$ is the minimum $\sigma_T(G,T)$ over all
spanning trees $T$ of $G$, and a spanning tree $T$ that minimizes
$\sigma_T(G,T)$ is called an {\it optimal tree}.

Let $P_n,C_n,K_n$ denote the path, the cycle, the complete graph,
respectively, on $n$ vertices. The {\it join} of two graphs $G$ and
$H$, denoted $G\vee H$, is the union of $G$ and $H$ and adding edges
from every vertex of $G$ to every vertex of $H$. For example,
$W_n=K_1\vee C_{n-1}$ is the wheel on $n$ vertices, $K_{m,n}=\bar
K_m\vee \bar K_n$ is the complete bipartite graph with $(m,n)$
partition. The {\it cartesian product} of two graphs $G$ and $H$,
denoted $G\times H$, is the graph with vertex set $V(G)\times V(H)$
and two vertices $(u,v)$ and $(u',v')$ are adjacent if and only if
either $[u=u'\, \mbox{and}\, vv'\in E(H)]$ or $[v=v'\, \mbox{and}\,
uu'\in E(G)]$. For example, $P_m\times P_n$ is the rectangular grid,
$C_m\times C_n$ is the torus grid.

A {\it block} of $G$ is a subgraph of $G$ which contains no cut
vertices and it is maximal with respect to this property. Two blocks
of $G$ have at most one vertex (a cut vertex) in common. As each
fundamental cycle is contained in a block, we have the following.

{\bf Proposition 2.1.} \ If $G$ has blocks $G_1,G_2,\ldots,G_k$,
then
$$\sigma_T(G)=\max_{1\leq i\leq k}\sigma_T(G_i).$$

So, we may assume that $G$ is itself a block, that is a 2-connected
graph (for $n\geq 3$). It is trivial that $\sigma_T(G)=1$ iff $G$ is
a tree. The {\it girth} of $G$ is the length of a shortest cycle in
$G$. By definition, we have a lower bound as follows.

{\bf Proposition 2.2.} \ Let $g(G)$ be the girth of $G$. Then
$\sigma_T(G)\geq g(G)-1$.

Several graphs attain this lower bound by choosing suitable spanning
trees. The following are some examples (see Figure 1, in which the
spanning trees are depicted by solid lines, while the cotrees by
dotted lines).

{\bf Proposition 2.3.} \ The following graphs have $\sigma_T(G)=g(G)-1$:\\
\indent (1) $\sigma_T(K_n)=2$ for the complete graphs $K_n$ ($n\geq 3$).\\
\indent (2) $\sigma_T(C_n)=n-1$ for the cycles $C_n$ ($n\geq 3$).\\
\indent (3) $\sigma_T(W_n)=2$ for the wheels $W_n=C_{n-1}\vee K_1$ ($n\geq 4$).\\
\indent (4) $\sigma_T(D_n)=2$ for the diamonds $D_n=K_2\vee \bar K_{n-2}$ ($n\geq 4$).\\
\indent (5) $\sigma_T(K_{m,n})=3$ for the complete bipartite graphs $K_{m,n}$ ($m,n\geq 2$).\\
\indent (6) $\sigma_T(P_3\times P_n)=3$ for special planar grids
$P_3\times P_n$ ($n\geq 2$).\\
\indent (7) $\sigma_T(G)=4$ for the Petersen graph $G$.

\begin{center}
\setlength{\unitlength}{0.4cm}
\begin{picture}(31,11)

\multiput(0,7)(1.5,0){2}{\circle*{0.3}}
\multiput(3.5,7)(1.5,0){2}{\circle*{0.3}}
\multiput(2.5,4)(0,6){2}{\circle*{0.3}}
\multiput(8,4)(0,3){3}{\circle*{0.3}}
\multiput(11,4)(0,3){3}{\circle*{0.3}}
\multiput(15,4)(2,0){4}{\circle*{0.3}}
\multiput(15,7)(2,0){4}{\circle*{0.3}}
\multiput(15,10)(2,0){4}{\circle*{0.3}}

\put(2.5,4){\line(-5,6){2.5}} \put(2.5,4){\line(-1,3){1}}
\put(2.5,4){\line(0,1){6}}\put(2.5,4){\line(1,3){1}}
\put(2.5,4){\line(5,6){2.5}} \put(8,4){\line(1,0){3}}
\put(8,4){\line(1,1){3}} \put(8,4){\line(1,2){3}}
\put(11,4){\line(-1,1){3}} \put(11,4){\line(-1,2){3}}
\put(15,4){\line(0,1){6}} \put(17,4){\line(0,1){6}}
\put(19,4){\line(0,1){6}} \put(21,4){\line(0,1){6}}
\put(15,7){\line(1,0){6}}

\bezier{12}(0,7)(1.25,8.5)(2.5,10) \bezier{12}(1.5,7)(2,8.5)(2.5,10)
\bezier{12}(3.5,7)(3,8.5)(2.5,10) \bezier{12}(5,7)(3.75,8.5)(2.5,10)
\bezier{12}(8,7)(9.5,7)(11,7) \bezier{12}(8,10)(9.5,8.5)(11,7)
\bezier{12}(8,7)(9.5,8.5)(11,10) \bezier{12}(8,10)(9.5,10)(11,10)
\bezier{30}(15,4)(18,4)(21,4) \bezier{30}(15,10)(18,10)(21,10)

\multiput(25.3,4)(4.4,0){2}{\circle*{0.3}}
\multiput(26.1,5.2)(2.8,0){2}{\circle*{0.3}}
\multiput(25.4,7.5)(4.2,0){2}{\circle*{0.3}}
\multiput(24,8)(7,0){2}{\circle*{0.3}}
\multiput(27.5,9)(1,0){1}{\circle*{0.3}}
\multiput(27.5,10.5)(1,0){1}{\circle*{0.3}}

\put(27.5,9){\line(0,1){1.5}} \qbezier(24,8)(25.75,9.25)(27.5,10.5)
\qbezier(27.5,10.5)(29.25,9.25)(31,8)
\qbezier(26.1,5.2)(26.8,7.1)(27.5,9)
\qbezier(28.9,5.2)(28.2,7.1)(27.5,9)
\qbezier(24,8)(24.7,7.75)(25.4,7.5)
\qbezier(29.6,7.5)(30.3,7.75)(31,8)
\qbezier(25.3,4)(25.7,4.6)(26.1,5.2)
\qbezier(28.9,5.2)(29.3,4.6)(29.7,4)

\bezier{15}(25.3,4)(27.5,4)(29.7,4)
\bezier{15}(25.4,7.5)(28,7.5)(29.6,7.5)
\bezier{15}(25.4,7.5)(27.15,6.35)(28.9,5.2)
\bezier{15}(26.1,5.2)(27.85,6.35)(29.6,7.5)
\bezier{15}(25.3,4)(24.65,6)(24,8)
\bezier{15}(31,8)(30.35,6)(29.7,4)

\put(-0.5,2.2){\makebox(1,0.5)[l]{\small (a) Diamond $D_6$}}
\put(8,2.2){\makebox(1,0.5)[l]{\small (b) $K_{3,3}$}}
\put(15.5,2.2){\makebox(1,0.5)[l]{\small (c) $P_3\times P_4$}}
\put(23.5,2.2){\makebox(1,0.5)[l]{\small (d) Petersen graph}}
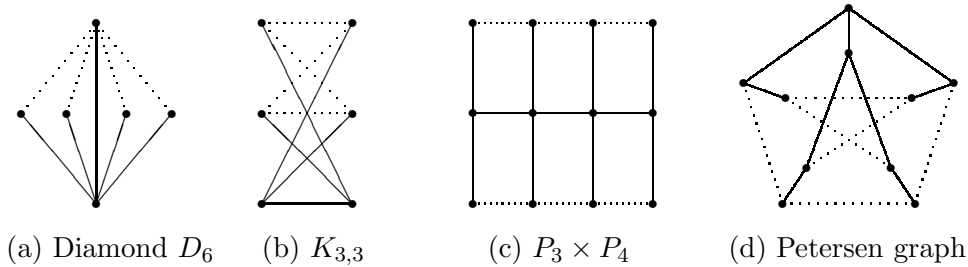
\put(8,0.5){\makebox(1,0.5)[l]{\small Figure 1. Examples in
Proposition 2.3.}}
\end{picture}
\end{center}

{\bf Proof.} \ (1) The complete graph $K_n$ ($n\geq 3$) has girth
$g(K_n)=3$ and a star $K_{1,n-1}$ is an optimal tree. (2) The cycle
$C_n$ ($n\geq 3$) has the unique fundamental cycle itself. (3) The
wheel $W_n=C_{n-1}\vee K_1$ has girth $3$ and the star $K_{1,n-1}$
is an optimal tree. (4) The diamond $D_n$ has girth $3$ and the star
$K_{1,n-1}$ is an optimal tree (see Figure 1(a)). (5) Let $G$ be a
complete bipartite graph $K_{m,n}$ with bipartition $(X,Y)$ where
$|X|=m,|Y|=n$ ($m,n\geq 2$). Then $G$ has girth $4$. We can
construct a spanning tree $T$ by taking a star $K_{1,n}$ with center
$x\in X$ and a star $K_{1,m}$ with center $y\in Y$ (which is called
a {\it double star} with diameter three, see Figure 1(b)). Then each
fundamental cycle with respect to $T$ has length $4$, and thus $T$
is optimal. (6) For the planar grid $P_3\times P_n$, the girth is 4
and the `caterpillar' with leaves on the boundary of outer face is
an optimal tree (see Figure 1(c)). (7) For the Petersen graph $G$,
the girth is 5 and we take the spanning tree $T$ as shown in Figure
1(d). Then every fundamental cycle with respect to $T$ has length 5.
This completes the proof. $\Box$

It is interesting to characterize the graphs satisfying Proposition
2.3, namely, those graphs having a spanning tree that every
fundamental cycle is a shortest cycle. We shall see more examples in
the next section.

\section{Characterization of low stretch graphs}
\hspace*{0.5cm} This section is intended to approach the open
problem of characterizing $\sigma_T(G)=3$. Madanlel et al.
\cite{Madan96} showed that $\sigma_T(G)\leq 3$ for all interval and
permutation graphs, and that a regular bipartite graph $G$ has
$\sigma_T(G)\leq 3$ if and only if it is complete. Moreover,
Brandst\"adt et al. \cite{Brand07} showed $\sigma_T(G)=3$ for
bipartite ATE-free graphs and convex graphs. Here, an ATE
(asteroidal triple of edges) in a graph $G$ is a set $A$ of three
edges that for any two edges $e_1,e_2\in A$, there is a path from
$e_1$ to $e_2$ that avoids the neighborhood of the third edge $e_3$
(the neighborhood of $uv$ is $N_G(u)\cup N_G(v)$). An ATE-free
(asteroidal-triple-edge-free) graph is one which does not contain
any ATE. The bipartite convex graphs form a special class of
bipartite ATE-free graphs. A bipartite graph $G$ with bipartition
$(X,Y)$ is said to be {\it convex} if $Y$ can be ordered as
$Y=\{y_1,y_2,\ldots,y_n\}$ such that the neighbor set $N_G(x_i)$ is
a consecutive sequence in $Y$ for each $x_i\in X$. We present more
results in this context.

\subsection{Complete $k$-partite graphs}
\hspace*{0.5cm} Let $\{V_1,V_2,\ldots, V_k\}$ be a partition of
$V(G)$ with $n_i=|V_i|$ ($1\leq i\leq k$). The complete $k$-partite
graph $K_{n_1,n_2,\ldots,n_k}$ with $k\geq 2$ is a graph such that
$uv\in E(G)$ if and only if $u\in V_i$ and $v\in V_j$ for $i\neq j$.

{\bf Theorem 3.1.} \ Suppose that $n_1\leq n_2\leq\cdots \leq n_k$
and $k\geq 3$. Then
$$\sigma_T(K_{n_1,n_2,\ldots,n_k})=\begin{cases}
2,& \mbox{if}\,\, n_1=1\\
3,& \mbox{otherwise}.
\end{cases}$$

{\bf Proof.} \ Let $G=K_{n_1,n_2,\ldots,n_k}$ ($k\geq 3$).
Obviously, the girth of $G$ is 3. When $n_1=1$, we can construct a
spanning tree as a star centered at the unique vertex of $V_1$. Then
all fundamental cycles are triangles, and thus $\sigma_T(G)=2$. When
$n_1\geq 2$, we will show that for any spanning tree $T$,
$\sigma_T(G,T)\geq 3$. By letting $X=V_2\cup\cdots \cup V_k$, we
have a complete bipartite graph $G'$ of bipartition $(V_1,X)$. There
are two cases to consider.\\
{\indent}(i) The spanning tree $T$ contains no edges between
vertices in $X$. Then $T$ is a spanning tree of $G'$ and a
fundamental cycle with respect to $T$ in $G'$ is one in $G$. As $G'$
is bipartite, a fundamental cycle in $G'$ has length at least
$4$, whence $\sigma_T(G,T)\geq 3$.\\
{\indent}(ii) The spanning tree $T$ contains some edges between
vertices in $X$. Suppose that $xy\in T$ with $x\in V_i$ and $y\in
V_j$ ($2\leq i<j\leq k$). Let $u\in V_1$ be such that
$d_T(u,x)<d_T(u,y)$. If $d_T(u,x)\geq 2$, then $d_T(u,y)\geq 3$,
thus $\sigma_T(G,T)\geq 3$. Otherwise $ux\in T$. Take $z\in
V_i,z\neq x$. Then $d_T(x,z)\geq 2$. If the path $P_{xz}$ in $T$
contains $u$, then $d_T(y,z)\geq 3$. Otherwise $d_T(u,z)\geq 3$,
whence $\sigma_T(G,T)\geq 3$.\\
{\indent} On the other hand, we can construct a spanning tree $T$ in
the complete bipartite graph $G'$ as a double star (as in
Proposition 2.3(5)). Then for an edge between the vertices of $V_1$
and $X$, the fundamental cycle has length four, while for an edge
between the vertices of $X$, the fundamental cycle has length three.
Thus $\sigma_T(G,T)=3$. This completes the proof. $\Box$

\subsection{Split graphs}
\hspace*{0.5cm}A graph $G$ is a {\it split graph} if its vertex set
$V(G)$ can be partitioned into a clique $X$ of $G$ and an
independent set $Y$ of $G$. For split graphs, \cite{Okamoto11}
showed that the spanning tree congestion problem is NP-complete.
However, the dual problem is easy. It has been known in
\cite{Brand04,Venka97} that $\sigma_T(G)\leq 3$ for split graphs
$G$. Here we describe a precise characterization as follows.

{\bf Theorem 3.2.} \ For a split graph $G$ (apart from a tree),
$\sigma_T(G)=2$ if and only if there exists a vertex $x_0\in X$ such
that every vertex $y\in Y\setminus N_G(x_0)$ is a pendant vertex (of
degree one). Otherwise $\sigma_T(G)=3$.

{\bf Proof.} \ If there exists a vertex $x_0\in X$ such that every
vertex $y\in Y\setminus N_G(x_0)$ is pendant, then we can construct
a spanning tree $T^*$ by the star with edges from $x_0$ to
$N_G(x_0)$, and by joining each remaining vertex $y\in Y\setminus
N_G(x_0)$ to its unique neighbor in $X$. Then for any $x,x'\in X$,
we have $d_{T^*}(x,x')\leq 2$. For any edge $xy\in E(G)$ with $x\in
X$ and $y\in N_G(x_0)$, the path between $x$ and $y$ in $T^*$ is
either $x_0y$ or $xx_0y$, thus $d_{T^*}(x,y)\leq 2$. For any edge
$xy\in E(G)$ with $x\in X$ and $y\in Y\setminus N_G(x_0)$, we have
$x\neq x_0$. Then $x$ is the unique neighbor of $y$, thus
$d_{T^*}(x,y)=1$. Therefore $\sigma_T(G,T^*)=2$ and so
$\sigma_T(G)=2$.

Conversely, if $\sigma_T(G)=2$, then there is a spanning tree $T$
such that $\sigma_T(G,T)=2$. This spanning tree $T$ restricted in
$G[X]$ must be a star with center $x_0$. For otherwise there would
be $x,x'\in X$ such that $d_T(x,x')\geq 3$. If a vertex $y\in Y
\setminus N_G(x_0)$ is adjacent to two vertices $x_1,x_2\in X$
(where $yx_1\in T$), then the fundamental cycle $yx_1x_0x_2$ has
length greater than three, which contradicts $\sigma_T(G,T)=2$.

Furthermore, we show that $\sigma_T(G)\leq 3$ in any case. To this
end, we construct a spanning tree $T$ as follows. We choose a vertex
$x_0\in X$ arbitrarily and take the star from $x_0$ to $N_G(x_0)$,
and join each vertex $y\in Y\setminus N_G(x_0)$ to a neighbor in
$X$. For any $x,x'\in X$, we have $d_{T}(x,x')\leq 2$. For any edge
$xy\in \bar T$ with $x\in X$ and $y\in N_G(x_0)$, the path between
$x$ and $y$ in $T$ is $xx_0y$, thus $d_{T}(x,y)=2$. If there is an
edge $xy\in \bar T$ with $x\in X$ and $y\in Y\setminus N_G(x_0)$,
and $yx'\in T$, then the path between $x$ and $y$ in $T$ is
$xx_0x'y$. Thus $d_{T}(x,y)=3$. Therefore, $\sigma_T(G,T)\leq 3$ and
so $\sigma_T(G)\leq 3$. This completes the proof. $\Box$

\subsection{Generalized convex graphs}
\hspace*{0.5cm} A bipartite graph $G$ with bipartition $(X,Y)$ is a
{\it chain graph} if there is an order $x_1,x_2,\ldots,x_m$ in $X$
such that $N_G(x_1)\subseteq N_G(x_2) \subseteq \cdots \subseteq
N_G(x_m)$. Previously, \cite{Okamoto11} showed that the minimum
congestion spanning tree problem is NP-hard even for chain graphs.
However, the counterpart in the tree-stretch problem is quite easy,
since a chain graph is a special convex graph and $\sigma_T(G)\leq
3$ in known in $\cite{Brand07}$.

Now we consider a generalization of convex graphs. A subset family
${\cal F}$ is called {\it laminar} (or {\it nested}) if for any two
sets $A,B\in {\cal F}$, at least one of $A\setminus B,B\setminus A,
A\cap B$ is empty, that is, $A\cap B\neq \emptyset \Rightarrow
A\subseteq B \, \mbox{or}\, B\subseteq A$.

{\bf Definition.} \ A bipartite graph $G$ with bipartition $(X,Y)$
is a {\it generalized convex graph} if there exists a tree $\tau
(Y)$ on the vertex set $Y$ such that for each $x_i\in X$, the
neighbor set $Y_i=N_G(x_i)$ induces a subpath in $\tau (Y)$ and the
subset family $\Sigma=\{Y_i: x_i\in X\}$ satisfies the following \\
{\indent} {\bf Laminar property} \ For each maximal subset $Y_0\in
\Sigma$ (there exists no $Y_i\in \Sigma$ such that $Y_0\subset
Y_i$), the subset family $\{Y_i\setminus Y_0: Y_i\cap Y_0 \neq
\emptyset,Y_i\in \Sigma\}$ is laminar.

For a convex graph $G$, $\tau (Y)$ is itself a path and the subset
family $\Sigma=\{Y_i: x_i\in X\}$ can be regarded as a set of
intervals on the line of $\tau (Y)$. For each maximal interval
$Y_0\in \Sigma$, if $Y_i\cap Y_0 \neq \emptyset, Y_j\cap Y_0 \neq
\emptyset$, then $Y_i\setminus Y_0$ and $Y_j\setminus Y_0$ is either
disjointed or one is included in another. Hence the subset family
$\{Y_i\setminus Y_0: Y_i\cap Y_0 \neq \emptyset,Y_i\in \Sigma\}$ is
laminar. Thus the above definition is indeed a generalization of
that of convex graphs. Moreover, a generalized convex graph is not
necessarily an ATE-free graph. For example, when $\tau (Y)$ is not a
path, let $y_1,y_2,y_3$ be three leaves (pendant vertices) of $\tau
(Y)$ in different branches such that there is a path from $y_i$ to
$y_j$ that avoids the neighborhood of $y_k$ (for $\{i,j,k\}=
\{1,2,3\}$). Then the three edges $e_1,e_2,e_3$ incident with
$y_1,y_2,y_3$, respectively, in $G$ constitute an ATE.

We are going to show that $\sigma_T(G)=3$ for generalized convex
graphs. Since a bipartite graph (apart from a tree) has girth
$g(G)=4$, we have $\sigma_T(G)\geq 3$. It suffices to construct an
optimal spanning tree with stretch three.

Let $\Sigma=\{Y_1,Y_2,\ldots,Y_m\}$ be the family of neighbor sets,
where $Y_i=N_G(x_i)$ for $x_i\in X$ ($1\leq i\leq m$). By
assumption, we are given a tree $\tau (Y)$ on $Y$ that each $Y_i$
induces a subpath of it ($1\leq i\leq m$). Suppose that $Y_1$
contains a leaf (pendant vertex) of $\tau (Y)$ and it is maximal in
$\Sigma$ in the sense of inclusion. We consider this leaf as the
root of the tree. Starting with $Y_1$, we define the {\it level
sets} $L_k$ in $\Sigma$ by the following procedure:\\
{\indent} (i) Define $L_1:=\{Y_1\}$. Set $\Sigma:=\Sigma\setminus
L_1$ and $k:=1$. \\
{\indent} (ii) For each $Y_i\in L_k$, if $Y_j\in \Sigma$ satisfies
that $Y_i\cap Y_j\neq \emptyset$, $Y_j\setminus Y_i\neq \emptyset$,
and $Y_i\cup Y_j$ is maximal (i.e., there is no other $Y_l$ such
that $Y_i\cap Y_l\neq \emptyset$, $Y_l\setminus Y_i\neq \emptyset$,
and $Y_i\cup Y_j \subset Y_i\cup Y_l$), then $Y_j$ is called a {\it
successor} is $Y_i$ (and $Y_i$ is the {\it predecessor} of $Y_j$).
Let $L_{k+1}$ be the set of successors of $Y_i$ for all $Y_i\in L_k$. \\
{\indent} (iii) Set $\Sigma:=\Sigma\setminus L_{k+1}$ and $k:=k+1$.
If $\bigcup_{Y_i\in L_1\cup L_2\cup \cdots\cup L_k}Y_i=Y$, then let
$h:=k$ and stop, else go to (ii).

By this procedure, we construct the level sets $L_1,L_2,\ldots,L_h$.
Let $\Sigma^*:=\bigcup_{1\leq k\leq h}L_k$, which is a subfamily of
$\Sigma$. For all neighbor sets $Y_i$ in $\Sigma^*$, no one is
contained in another, and they constitute a cover of $Y$. Also, they
can be regarded as a directed tree rooted at $Y_1$ and running down
level by level. If $Y_j$ and $Y_l$ are successors of $Y_i$ in this
directed tree, then by the laminar property, we see that
$(Y_j\setminus Y_i)\cap (Y_l\setminus Y_i)=\emptyset$. Also, for
$Y_i\in L_{k-1}$ and $Y_j\in L_{k+1}$, we have $Y_i\cap
Y_j=\emptyset$.

In this situation, there may be some neighbor sets $Y_q\in \Sigma
\setminus \Sigma^*$, which are discarded in the above procedure. For
each $Y_q\in \Sigma \setminus \Sigma^*$, there must be a $Y_i\in
L_k$ and its successor $Y_j\in L_{k+1}$ such that $Y_i\cap Y_q\neq
\emptyset$, and $Y_i\cup Y_q \subseteq Y_i\cup Y_j$. For otherwise
we may choose $Y_q$ in the above procedure.

By means of the level structure $\{L_1,L_2,\ldots,L_h\}$, we
construct the spanning tree $T$ by the following algorithm.

\vskip 3mm {\bf Construction Algorithm}\\
\hspace*{0.2cm}\line(1,0){170}\\
{\indent} (1) For $L_1=\{Y_1\}$, construct a star $T_1$ with center
$x_1$ and all leaves $y\in Y_1$. Set $T:=T_1$ and $k:=1$.\\
{\indent} (2) For each neighbor set $Y_i\in L_k$, consider a
successor $Y_j\in L_{k+1}$, and construct a star $T_j$ with center
$x_j$ and all leaves $y\in (Y_j\setminus Y_i)\cup \{\bar y\}$, where
$\bar y$ is the last vertex in $Y_i\cap Y_j$ (according to the order
of the path of $Y_i$). Set $T:=T\cup T_j$. Repeat this step for all
successors of $Y_i$ and all $Y_i\in L_k$. \\
{\indent} (3) Set $k:=k+1$. If $k<h$, then go to (2). \\
{\indent} (4) For each neighbor set $Y_q\in \Sigma\setminus
\Sigma^*$, suppose that $Y_i\cap Y_q\neq \emptyset$ and
$Y_q\subseteq Y_i\cup Y_j$ for some $Y_i\in L_k$ and its successor
$Y_j\in L_{k+1}$. Let $\bar y$ be the last vertex in $Y_i\cap Y_q$.
Then set $T:=T\cup \{x_q\bar y\}$.\\
\hspace*{0.2cm}\line(1,0){170}\\

We claim that the output $T$ of the above algorithm is indeed a
spanning tree of $G$. In fact, we first construct a star $T_1$ with
center $x_1$ for level $L_1$. When $Y_i\in L_k$ has been considered,
we have a star $T_i$ with center $x_i$. Then we consider a successor
$Y_j\in L_{k+1}$ of $Y_i$ and add a star $T_j$. Since the stars
$T_i$ and $T_j$ have only one leaf in common, $T_i\cup T_j$ is
connected and contains no cycles, and thus is a tree. If $Y_i$ has
another successor $Y_l$, then by the laminar property, we have
$(Y_j\setminus Y_i)\cap (Y_l\setminus Y_i)=\emptyset$. Then $T_l$
and $T_i$ have one leaf in common, $T_l$ and $T_j$ have at most one
leaf in common (if the $\bar y\in Y_i$ is the same for $T_l$ and
$T_j$). Hence $T_i\cup T_j\cup T_l$ is also a tree. In this way, we
construct a set of stars in which any two stars have at most one
leaf in common. So we obtain a tree $T$ in Steps (1)-(3). In Step
(4), we add more pendant edges (with new leaves $x_q$) to $T$.
Additionally, all vertices of $G$ are considered when the algorithm
terminates. Therefore $T$ is finally a spanning tree.

{\bf Theorem 3.3.} \ For a generalized convex graph $G$ (apart from
a tree), it holds that $\sigma_T(G)=3$.

{\bf Proof.} \ We proceed to show that the spanning tree $T$
constructed by the above algorithm has stretch three. For each
cotree-edge $e\in \bar T$, there are two cases to consider:

{\bf Case 1:} \ $e=x_jy$ with $Y_j\in L_k$ for some $L_k$ in
$\sigma^*$. Let $Y_i\in L_{k-1}$ be the predecessor of $Y_j$. Then
$y\in Y_i\cap Y_j$. Thus $x_iy,x_i\bar y$ and $x_j\bar y$ are
contained in $T$ (where $\bar y$ is the last vertex in $Y_i\cap
Y_j$). Hence $e=x_jy$ and these three edges in $T$ constitute the
fundamental cycle with respect to $e$, which has length four.

{\bf Case 2:} \ $e=x_qy$ with $Y_q\in \Sigma\setminus \Sigma^*$.
Then there is some $Y_i\in L_k$ and its successor $Y_j\in L_{k+1}$
such that $Y_i\cap Y_q\neq \emptyset$ and $Y_q\subseteq Y_i\cup
Y_j$. If $y\in Y_i$ (say $Y_q\subseteq Y_i$), then $x_iy,x_i\bar y,
x_q\bar y\in T$ (where $\bar y$ is the last vertex in $Y_i\cap
Y_q$). Thus $e=x_qy$ and these three edges in $T$ constitute the
fundamental cycle with respect to $e$, which has length four. If
$y\in Y_j\setminus Y_i$, then by the laminar property, $Y_q\setminus
Y_i\subseteq Y_j\setminus Y_i$. Let $\bar y$ be the last vertex in
$Y_i\cap Y_j$. Then $\bar y\in Y_q$ and $x_jy,x_j\bar y,x_q\bar y\in
T$. Thus these three edges in $T$ and $e=x_qy\notin T$ also yield a
length four fundamental cycle.

To summarize, for every cotree-edge $e\in \bar T$, the fundamental
cycle with respect to $e$ has length four. Therefore,
$\sigma_T(G,T)=3$ and the theorem is proved. $\Box$

\section{Planar grids}
\hspace*{0.5cm} It is known that the minimum stretch spanning tree
problem is NP-hard for planar graphs in general \cite{Fekete01}. We
discuss some planar grids in this section.

Let $G$ be a simple connected planar graph. Suppose that we have a
planar embedding of $G$ on the plane so that it is a {\it plane
graph}. For a face $f$ of $G$, the degree of $f$, denoted by $d(f)$,
is the number of edges in its boundary. Our approach is based on the
spanning trees of the dual graph. The dual graph $G^*$ of $G$ is
defined as follows. Each face $f$ of $G$ (including the outer face)
corresponds to a vertex $f^*$ in $G^*$, and each edge $e$ of $G$
corresponds to an edge $e^*$ of $G^*$ in such a way that two
vertices $f^*$ and $g^*$ are joined by an edge $e^*$ in $G^*$ if and
only if their corresponding faces $f$ and $g$ are separated by the
edge $e$ in $G$. We may place each vertex $f^*$ in the face $f$ of
$G$ and draw each edge $e^*$ to cross the edge $e$ of $G$ exactly
once. This dual graph $G^*$ is also a plane graph.

A prominent property of duality is: A cycle $C$ of $G$ corresponds
an edge-cut (cocycle) $C^*$ of $G^*$, and an edge-cut $B$ of $G$
corresponds a cycle $B^*$ of $G^*$. In particular, for a spanning
tree $T$ of $G$, the cotree $\bar T$ corresponds to a spanning tree
$\bar T^*$ of $G^*$. A fundamental cycle with respect to $T$ in $G$
corresponds to a fundamental edge-cut with respect to $\bar T^*$ in
$G^*$ (see \cite{Bondy08} for details). For example, the cube $Q_3$
is shown in Figure 2(a) and a spanning tree $T$ with solid lines in
Figure 2(b). Meanwhile, the spanning tree $\bar T^*$ with dotted
lines of the dual graph $G^*$ is also drawn in Figure 2(b), in which
the vertices of faces are represented by small circles and the
vertex of outer face is denoted by $O$.

\begin{center}
\setlength{\unitlength}{0.3cm}
\begin{picture}(32,16)

\multiput(1,4)(10,0){2}{\circle*{0.3}}
\multiput(4,7)(4,0){2}{\circle*{0.3}}
\multiput(4,11)(4,0){2}{\circle*{0.3}}
\multiput(1,14)(10,0){2}{\circle*{0.3}}
\multiput(16,4)(10,0){2}{\circle*{0.3}}
\multiput(19,7)(4,0){2}{\circle*{0.3}}
\multiput(19,11)(4,0){2}{\circle*{0.3}}
\multiput(16,14)(10,0){2}{\circle*{0.3}}
\multiput(17.5,9)(7,0){2}{\circle{0.3}}
\multiput(21,5.5)(0,3.5){3}{\circle{0.3}}
\multiput(30,9)(3,0){1}{\circle{0.3}}

\put(1,4){\line(1,0){10}} \put(4,7){\line(1,0){4}}
\put(4,11){\line(1,0){4}}\put(1,14){\line(1,0){10}}
\put(1,4){\line(0,1){10}} \put(4,7){\line(0,1){4}}
\put(8,7){\line(0,1){4}} \put(11,4){\line(0,1){10}}
\put(1,4){\line(1,1){3}} \put(8,11){\line(1,1){3}}
\put(1,14){\line(1,-1){3}} \put(8,7){\line(1,-1){3}}
\put(19,7){\line(1,0){4}} \put(19,11){\line(1,0){4}}
\put(19,7){\line(0,1){4}} \put(16,4){\line(1,1){3}}
\put(23,11){\line(1,1){3}} \put(16,14){\line(1,-1){3}}
\put(23,7){\line(1,-1){3}}

\bezier{30}(21.2,9)(25,9)(29.8,9)
\bezier{40}(21,12.6)(24,20)(30,9.1)
\bezier{40}(21,5.4)(25,-2)(30,8.9) \bezier{25}(17.4,9)(13,11)(15,14)
\bezier{55}(15,14)(25,21)(30,9.1)

\put(5.6,5.2){\makebox(1,0.5)[l]{\footnotesize $1$}}
\put(2.2,8.5){\makebox(1,0.5)[l]{\footnotesize $1$}}
\put(5.6,8.5){\makebox(1,0.5)[l]{\footnotesize $2$}}
\put(9.2,8.5){\makebox(1,0.5)[l]{\footnotesize $1$}}
\put(5.6,12.3){\makebox(1,0.5)[l]{\footnotesize $1$}}
\put(31,8.5){\makebox(1,0.5)[l]{\footnotesize $O$}}
\put(1.5,0.8){\makebox(1,0.5)[l]{\small (a) The cube $Q_3$}}
\put(13,0.8){\makebox(1,0.5)[l]{\small (b) Spanning trees of $Q_3$
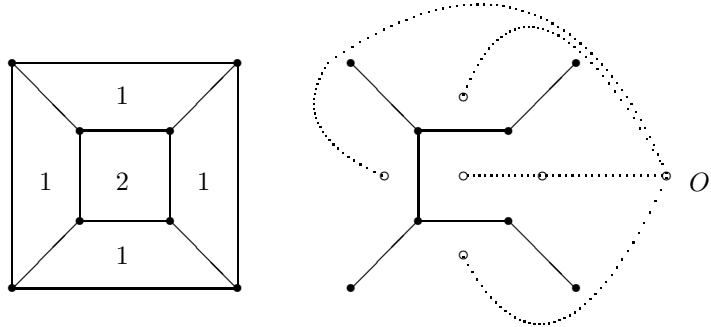
and its dual}} \put(5,-1){\makebox(1,0.5)[l]{\small Figure 2. The
cube $Q_3$ and its spanning trees.}}
\end{picture}
\end{center}

For a face $f$ of plane graph $G$ (a vertex of $G^*$), we define the
{\it level} of $f$, denoted by $\lambda(f)$, to be the length of a
shortest path from the vertex $f$ to the vertex $O$ of outer face in
$G^*$. We denote by $L_i$ the set of faces having level $i$ ($i=0,
1,\ldots$). Then the levels can be determined by the following
procedure:\\
{\indent} (i) Let $\lambda(O)=0$ and $L_0=\{O\}$.\\
{\indent} (ii) If $L_i$ has been defined, then for any face $f$
whose level $\lambda(f)$ is not defined and it is \\
{\indent} $\quad$ adjacent to a face $g\in L_i$, set $\lambda(f)=i+1$.\\
For example, the levels of the faces in $Q_3$ are shown in Figure
2(a) by the number in each face (except $O$ with level $0$), where
$|L_1|=4,|L_2| =1$. Here, we first consider the outer vertex $O$ as
the root. Then, all vertices in $L_1$ have the same {\it
predecessor} $O$. In general, when $\lambda(f)=i+1$ is defined in
terms of an adjacent vertex $g\in L_i$, $g$ is the predecessor of
$f$. Thus, a rooted tree (called {\it search tree}) is obtained
level by level. In this respect, we define the {\it maximum level}
of $G$ by
$$\lambda_{\max}(G):=\max_{f\in F}\lambda(f),$$ where $F$ is the set
of the faces of $G$. This is the height of the search tree.

\subsection{Rectangular grids}
\hspace*{0.5cm} First, we consider the rectangular grids $G=
P_m\times P_n$ ($2\leq m\leq n$) on the plane. Let $V(G):=
\{(i,j):1\leq i\leq m,1\leq j\leq n\}$ denote the vertex set of $G$,
and $(i,j)$ is adjacent to $(i',j')$ if $|i-i'|+|j-j|=1$ (see Figure
3(a)). Similar to the notation of matrices, we may call $R_i:=
\{(i,j):1\leq j\leq n\}$ the $i$-th row, and $Q_j:=\{(i,j): 1\leq
i\leq m\}$ the $j$-th column. The edges in the rows are called {\it
horizontal edges}. The edges in the columns are called {\it vertical
edges}.

\begin{center}
\setlength{\unitlength}{0.32cm}
\begin{picture}(36,15.5)

\multiput(1,4)(3,0){5}{\circle*{0.3}}
\multiput(1,7)(3,0){5}{\circle*{0.3}}
\multiput(1,10)(3,0){5}{\circle*{0.3}}
\multiput(1,13)(3,0){5}{\circle*{0.3}}
\multiput(19,4)(3,0){5}{\circle*{0.3}}
\multiput(19,7)(3,0){5}{\circle*{0.3}}
\multiput(19,10)(3,0){5}{\circle*{0.3}}
\multiput(19,13)(3,0){5}{\circle*{0.3}}
\multiput(20.5,5.5)(3,0){4}{\circle{0.3}}
\multiput(20.5,8.5)(3,0){4}{\circle{0.3}}
\multiput(20.5,11.5)(3,0){4}{\circle{0.3}}
\multiput(36,8.5)(3,0){1}{\circle{0.3}}

\put(1,4){\line(1,0){12}} \put(1,7){\line(1,0){12}}
\put(1,10){\line(1,0){12}}\put(1,13){\line(1,0){12}}
\put(1,4){\line(0,1){9}} \put(4,4){\line(0,1){9}}
\put(7,4){\line(0,1){9}} \put(10,4){\line(0,1){9}}
\put(13,4){\line(0,1){9}}

\put(19,7){\line(1,0){12}} \put(19,4){\line(0,1){3}}
\put(19,10){\line(0,1){3}} \put(19,10){\line(1,0){3}}
\put(22,4){\line(0,1){9}} \put(25,4){\line(0,1){9}}
\put(28,4){\line(0,1){3}} \put(28,10){\line(1,0){3}}
\put(28,10){\line(0,1){3}} \put(25,10){\line(1,0){3}}
\put(28,13){\line(1,0){3}} \put(28,4){\line(1,0){3}}

\bezier{30}(20.5,8.5)(16,9.5)(18,13)
\bezier{50}(18,13)(21,18)(31,15.5)
\bezier{40}(31,15.5)(34,15)(36,8.5)

\bezier{13}(23.5,8.5)(23.5,10.2)(23.5,11.5)
\bezier{13}(26.5,8.5)(28,8.5)(29.5,8.5)
\bezier{25}(29.5,8.5)(32,8.5)(36,8.5)
\bezier{50}(20.5,5.5)(21,1.5)(31,2)
\bezier{43}(23.5,5.5)(24,2)(31,2.5)
\bezier{36}(26.5,5.5)(27,2.5)(31,3)
\bezier{28}(29.5,5.5)(32,5.5)(36,8.5)
\bezier{50}(20.5,11.5)(21,17)(31,15)
\bezier{43}(23.5,11.5)(24,16)(31,14.5)
\bezier{36}(26.5,11.5)(27,15)(31,14)
\bezier{28}(29.5,11.5)(32,11.5)(36,8.5)
\bezier{35}(31,1.98)(34.8,3.2)(36,8.5)
\bezier{34}(31,2.5)(34.4,4)(36,8.5)
\bezier{33}(31,3)(34,4.3)(36,8.5)
\bezier{35}(31,15)(34,14.4)(36,8.5)
\bezier{34}(31,14.5)(34,14)(36,8.5)
\bezier{33}(31,14)(33.5,13.5)(36,8.5)

\put(2.3,11.3){\makebox(1,0.5)[l]{\footnotesize $1$}}
\put(5.3,11.3){\makebox(1,0.5)[l]{\footnotesize $1$}}
\put(8.3,11.3){\makebox(1,0.5)[l]{\footnotesize $1$}}
\put(11.3,11.3){\makebox(1,0.5)[l]{\footnotesize $1$}}
\put(2.3,8.3){\makebox(1,0.5)[l]{\footnotesize $1$}}
\put(11.3,8.3){\makebox(1,0.5)[l]{\footnotesize $1$}}
\put(2.3,5.3){\makebox(1,0.5)[l]{\footnotesize $1$}}
\put(5.3,5.3){\makebox(1,0.5)[l]{\footnotesize $1$}}
\put(8.3,5.3){\makebox(1,0.5)[l]{\footnotesize $1$}}
\put(11.3,5.3){\makebox(1,0.5)[l]{\footnotesize $1$}}
\put(5.3,8.3){\makebox(1,0.5)[l]{\footnotesize $2$}}
\put(8.3,8.3){\makebox(1,0.5)[l]{\footnotesize $2$}}
\put(0,13.5){\makebox(1,0.5)[l]{\footnotesize $(1,1)$}}
\put(12.2,13.5){\makebox(1,0.5)[l]{\footnotesize $(1,5)$}}
\put(0,3.1){\makebox(1,0.5)[l]{\footnotesize $(4,1)$}}
\put(12.2,3.1){\makebox(1,0.5)[l]{\footnotesize $(4,5)$}}
\put(36.6,8.3){\makebox(1,0.5)[l]{\footnotesize $O$}}
\put(2,0.8){\makebox(1,0.5)[l]{\small (a) Grid $P_4\times P_5$}}
\put(15,0.8){\makebox(1,0.5)[l]{\small (b) Spanning trees of
$P_4\times P_5$ and its dual}} \put(7,-1){\makebox(1,0.5)[l]{\small
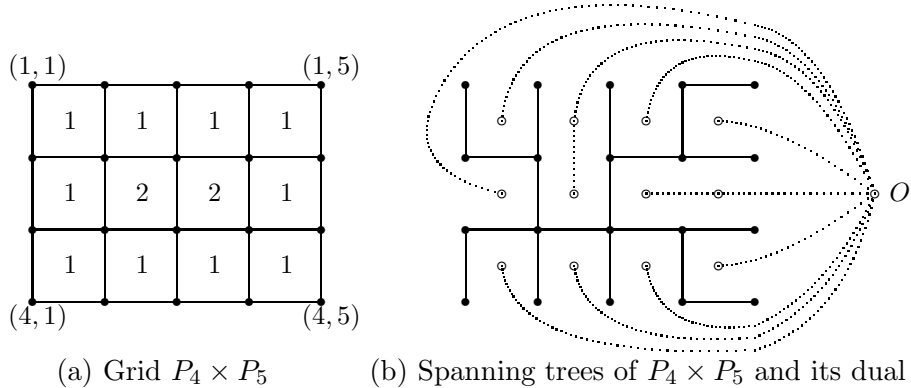
Figure 3. Grid $P_4\times P_5$ and spanning trees.}}
\end{picture}
\end{center}

Hruska \cite{Hruska08} proved the tree-congestion as follows ($m\leq
n$):
$$c_T(P_m\times P_n)=\begin{cases}
m, & \mbox{if } m=n \mbox{ or } m \mbox{ odd}\\
m+1, & \mbox{otherwise}.
\end{cases}$$
In the following we derive a similar formula for the tree-stretch:

{\bf Theorem 4.1.} \ For the rectangular grids $P_m\times P_n$ with
$2\leq m\leq n$, we have
$$\sigma_T(P_m\times P_n)=2\left\lfloor\frac{m}{2}\right\rfloor+1.$$

{\bf Proof.} \ Let $G=P_m\times P_n$ ($2\leq m\leq n$). We first
show that
$$\lambda_{\max}(G)=\left\lfloor\frac{m}{2}\right\rfloor.$$
By induction on $m$. When $m=2,3$, all faces have level $1$, so
$\lambda_{\max}(G)=1$ and the assertion holds. Assume that $m\geq 4$
and the assertion holds for smaller $m$. We delete the boundary of
the outer faces from $G$ (the vertices and the edges on this
boundary are deleted). Then the remaining graph is $G'=P_{m-2}\times
P_{n-2}$. In this transformation, all faces with level $1$ are
removed. Therefore $\lambda_{\max}(G)=\lambda_{\max}(G')+1$. By
induction hypothesis, $\lambda_{\max}(G')=\lfloor (m-2)/2\rfloor$.
Hence
$$\lambda_{\max}(G)=\left\lfloor\frac{m-2}{2}\right\rfloor+1=
\left\lfloor\frac{m}{2}\right\rfloor.$$ For example, the levels of
$P_4\times P_5$ are shown in Figure 3(a) and $\lambda_{\max}=\lfloor
m/2\rfloor=2$.

We next show the lower bound
\begin{equation}
\sigma_T(G)\geq 2\lambda_{\max}+1.
\end{equation}
In fact, let $T$ be any given spanning tree of $G$. Then the cotree
$\bar T$ determines a spanning tree $\bar T^*$ in $G^*$. Suppose
that $f_0$ is a face with the maximum level $\lambda_{\max}$. For
brevity, we still denote its vertex in $G^*$ by $f_0$ and write
$\lambda=\lambda_{\max}$. Then the distance between $f_0$ and $O$ in
$\bar T^*$ is at least $\lambda$. Let $P^*$ be the path from $f_0$
to $O$ in $\bar T^*$ with the last edge $e^*_0$ incident with $O$.
The tree-edge $e^*_0$ in the spanning tree $\bar T^*$ determines a
fundamental edge-cut $C^*=\partial(X_{e_0})$, where $\bar T^*-e^*_0$
has two components and $X_{e_0}$ is the vertex set of the component
containing $P^*$. Then this fundamental edge-cut $C^*$ with respect
to $\bar T^*$ in $G^*$ corresponds to a fundamental cycle $C$ with
respect to $T$ in $G$. So, this fundamental cycle $C$ is determined
by the cotree edge $e_0$ on the boundary of the outer face that
corresponds to the edge $e^*_0$ in $P^*$. Note that all faces in
$P^*$ (with labels $1,2,\ldots,\lambda$) are contained in the region
surrounded by $C$. Without loss of generality, assume that $e_0$ is
on the row $R_1$. We draw $\lambda$ horizontal straight lines
passing through the centers of square faces of $P^*$. Then each of
these straight lines intersects $C$ at two vertical edges. Besides,
$C$ must have at least two more horizontal edges. Hence $C$ has
length at least $2\lambda +2$. Consequently, for any spanning tree
$T$, we find a fundamental cycle $C$ with length at least $2\lambda
+2$. By the arbitrariness of $T$, the lower bound (3) is proved.

Conversely, we can construct a spanning tree $T^*$ by taking all
columns and the row $R_{\lfloor m/2\rfloor}$. Then the maximal
fundamental cycles have length $2\lfloor m/2\rfloor+2$. Thus the
spanning tree $T^*$ is optimal. This completes the proof. $\Box$

\subsection{Triangular grids}
\hspace*{0.5cm} We next consider the triangular grids $T_n$, which
is defined as follows. The vertex set can be represented as
$\{(x,y)\in {\mathbf Z}^2: x+y\leq n,x,y\geq 0\}$ on the plane, and
two vertices $(x,y)$ and $(x',y')$ are joined by an edge if
$|x-x'|+|y-y'|=1$ or $|x-x'|+|y-y'|=2$ and $x+y=x'+y'$ (refer to
\cite{Lin14}). For example, $T_4$ is shown in Figure 4, and $T_1$ (a
triangle), $T_2$ and $T_3$ are shown in Figure 5. In this plane
embedding of $T_n$, the straight-lines $\{(x,y)\in {\mathbf R}^2:
y=k\}$ ($0\leq k\leq n-1$) are called {\it horizontal lines} and the
edges on them are called {\it horizontal edges}. Symmetrically, the
straight-lines $\{(x,y)\in {\mathbf R}^2: x=k\}$ ($0\leq k\leq n-1$)
are called {\it vertical lines} and the edges on them are called
{\it vertical edges}. In addition, there are slant edges.

\begin{center}
\setlength{\unitlength}{0.45cm}
\begin{picture}(26,12)
\multiput(1,3)(2,0){5}{\circle*{0.3}}
\multiput(1,5)(2,0){4}{\circle*{0.3}}
\multiput(1,7)(2,0){3}{\circle*{0.3}}
\multiput(1,9)(2,0){2}{\circle*{0.3}}
\multiput(1,11)(2,0){1}{\circle*{0.3}}

\put(1,3){\line(1,0){8}} \put(1,5){\line(1,0){6}}
\put(1,7){\line(1,0){4}} \put(1,9){\line(1,0){2}}
\put(1,3){\line(0,1){8}} \put(3,3){\line(0,1){6}}
\put(5,3){\line(0,1){4}} \put(7,3){\line(0,1){2}}
\put(1,5){\line(1,-1){2}} \put(1,7){\line(1,-1){4}}
\put(1,9){\line(1,-1){6}} \put(1,11){\line(1,-1){8}}

\put(1.4,3.5){\makebox(1,0.5)[l]{\small $1$}}
\put(3.4,3.5){\makebox(1,0.5)[l]{\small $1$}}
\put(5.4,3.5){\makebox(1,0.5)[l]{\small $1$}}
\put(7.4,3.5){\makebox(1,0.5)[l]{\small $1$}}
\put(2.1,4.2){\makebox(1,0.5)[l]{\small $2$}}
\put(4.1,4.2){\makebox(1,0.5)[l]{\small $2$}}
\put(6.1,4.2){\makebox(1,0.5)[l]{\small $2$}}
\put(1.4,5.5){\makebox(1,0.5)[l]{\small $1$}}
\put(5.4,5.5){\makebox(1,0.5)[l]{\small $1$}}
\put(3.4,5.5){\makebox(1,0.5)[l]{\small $3$}}
\put(2.1,6.2){\makebox(1,0.5)[l]{\small $2$}}
\put(4.1,6.2){\makebox(1,0.5)[l]{\small $2$}}
\put(1.4,7.5){\makebox(1,0.5)[l]{\small $1$}}
\put(3.4,7.5){\makebox(1,0.5)[l]{\small $1$}}
\put(2.1,8.2){\makebox(1,0.5)[l]{\small $2$}}
\put(1.4,9.5){\makebox(1,0.5)[l]{\small $1$}}

\multiput(14,3)(2,0){5}{\circle*{0.3}}
\multiput(14,5)(2,0){4}{\circle*{0.3}}
\multiput(14,7)(2,0){3}{\circle*{0.3}}
\multiput(14,9)(2,0){2}{\circle*{0.3}}
\multiput(14,11)(2,0){1}{\circle*{0.3}}

\multiput(14.8,3.7)(2,0){4}{\circle{0.3}}
\multiput(15.4,4.4)(2,0){3}{\circle{0.3}}
\multiput(14.8,5.7)(2,0){3}{\circle{0.3}}
\multiput(15.4,6.4)(2,0){2}{\circle{0.3}}
\multiput(14.8,7.7)(2,0){2}{\circle{0.3}}
\multiput(15.4,8.4)(2,0){1}{\circle{0.3}}
\multiput(14.8,9.7)(2,0){1}{\circle{0.3}}
\multiput(21,10)(2,0){1}{\circle{0.3}}

\put(14,5){\line(1,0){6}} \put(14,7){\line(1,0){4}}
\put(14,9){\line(1,0){2}} \put(20,3){\line(1,0){2}}
\put(14,3){\line(0,1){2}} \put(14,9){\line(0,1){2}}
\put(16,3){\line(0,1){6}} \put(18,3){\line(0,1){2}}
\put(20,3){\line(0,1){2}}

\bezier{5}(14.9,3.8)(15,4)(15.3,4.3)
\bezier{5}(16.9,3.8)(17,4)(17.3,4.3)
\bezier{5}(18.9,3.8)(19,4)(19.3,4.3)
\bezier{5}(14.9,5.8)(15,6)(15.3,6.3)
\bezier{5}(16.9,5.8)(17,6)(17.3,6.3)
\bezier{8}(17.5,6.4)(18,6)(18.7,5.7)
\bezier{5}(14.9,7.8)(15,8)(15.3,8.3)
\bezier{40}(14.8,3.5)(16,1.5)(22,2)
\bezier{35}(16.8,3.5)(18,1.8)(22,2.3)
\bezier{30}(18.8,3.5)(19.6,2.1)(22,2.6)
\bezier{40}(14.7,5.7)(10.5,8.5)(13,11)
\bezier{30}(14.7,7.7)(11.8,9.5)(13.5,11)
\bezier{45}(22,2)(27.1,3.5)(21,10)
\bezier{40}(22,2.3)(26.1,4)(21,10)
\bezier{35}(22,2.6)(25,4.5)(21,10)
\bezier{45}(13,11)(15.5,13.3)(21,10)
\bezier{40}(13.5,11)(16,13)(21,10)
\bezier{40}(14.8,9.7)(17,9.8)(21,10)
\bezier{30}(16.8,7.7)(19,9)(21,10)
\bezier{30}(18.8,5.7)(20,7.5)(21,10)
\bezier{40}(20.8,3.7)(20.9,6.7)(21,10)

\put(1.4,3.5){\makebox(1,0.5)[l]{\small $1$}}
\put(3.4,3.5){\makebox(1,0.5)[l]{\small $1$}}
\put(5.4,3.5){\makebox(1,0.5)[l]{\small $1$}}
\put(7.4,3.5){\makebox(1,0.5)[l]{\small $1$}}
\put(2.1,4.2){\makebox(1,0.5)[l]{\small $2$}}
\put(4.1,4.2){\makebox(1,0.5)[l]{\small $2$}}
\put(6.1,4.2){\makebox(1,0.5)[l]{\small $2$}}
\put(1.4,5.5){\makebox(1,0.5)[l]{\small $1$}}
\put(5.4,5.5){\makebox(1,0.5)[l]{\small $1$}}
\put(3.4,5.5){\makebox(1,0.5)[l]{\small $3$}}
\put(2.1,6.2){\makebox(1,0.5)[l]{\small $2$}}
\put(4.1,6.2){\makebox(1,0.5)[l]{\small $2$}}
\put(1.4,7.5){\makebox(1,0.5)[l]{\small $1$}}
\put(3.4,7.5){\makebox(1,0.5)[l]{\small $1$}}
\put(2.1,8.2){\makebox(1,0.5)[l]{\small $2$}}
\put(1.4,9.5){\makebox(1,0.5)[l]{\small $1$}}
\put(21.5,9.8){\makebox(1,0.5)[l]{\small $O$}}
\put(1,0.8){\makebox(1,0.5)[l]{\small (a) Triangular grid $T_4$}}
\put(12,0.8){\makebox(1,0.5)[l]{\small (b) Spanning trees of $T_4$
and its dual}} \put(4,-1){\makebox(1,0.5)[l]{\small Figure 4.
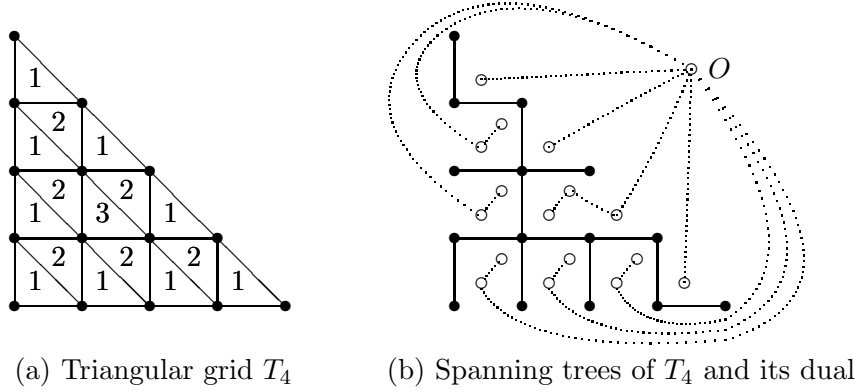
Triangular grid $T_4$ and spanning trees.}}
\end{picture}
\end{center}

Ostrovskii \cite{Ostrov10} developed an approach, called {\it
center-tail system}, to deal with the spanning tree congestion
problem for planar graphs, and obtained the result for triangular
grids as follows:
$$c_T(T_n)=\begin{cases}
4k, & \mbox{if } n=3k \\
4k, & \mbox{if } n=3k+1 \\
4k+2, & \mbox{if } n=3k+2
\end{cases}$$
(however, $T_n$ here is our $T_{n-1}$). We obtain the corresponding
result for tree-stretch as follows.

{\bf Theorem 4.2.} \ For the triangular grids $T_n$, we have
$$\sigma_T(T_n)=\left\lceil\frac{2n}{3}\right\rceil+1.$$

{\bf Proof.} \ For the triangular grids $T_n$, we first show that
$$\lambda_{\max}(T_n)=\left\lceil\frac{2n}{3}\right\rceil.$$
We use induction on $n$. When $1\leq n\leq 3$, the levels of faces
for $T_1,T_2,T_3$ are shown in Figure 5, in which
$\lambda_{\max}(T_1)=1,\lambda_{\max}(T_2)=\lambda_{\max}(T_3)=2$.
Hence the assertion holds for $1\leq n\leq 3$.

\begin{center}
\setlength{\unitlength}{0.45cm}
\begin{picture}(21,9)
\multiput(1,3)(3,0){2}{\circle*{0.3}}
\multiput(1,6)(2,0){1}{\circle*{0.3}}
\multiput(7,3)(2,0){3}{\circle*{0.3}}
\multiput(7,5)(2,0){2}{\circle*{0.3}}
\multiput(7,7)(2,0){1}{\circle*{0.3}}
\multiput(15,3)(2,0){4}{\circle*{0.3}}
\multiput(15,5)(2,0){3}{\circle*{0.3}}
\multiput(15,7)(2,0){2}{\circle*{0.3}}
\multiput(15,9)(2,0){1}{\circle*{0.3}}

\put(1,3){\line(1,0){3}} \put(1,3){\line(0,1){3}}
\put(1,6){\line(1,-1){3}}

\put(7,3){\line(1,0){4}} \put(7,5){\line(1,0){2}}
\put(7,3){\line(0,1){4}} \put(9,3){\line(0,1){2}}
\put(7,5){\line(1,-1){2}} \put(7,7){\line(1,-1){4}}

\put(15,3){\line(1,0){6}} \put(15,5){\line(1,0){4}}
\put(15,7){\line(1,0){2}} \put(15,3){\line(0,1){6}}
\put(17,3){\line(0,1){4}} \put(19,3){\line(0,1){2}}
\put(15,5){\line(1,-1){2}} \put(15,7){\line(1,-1){4}}
\put(15,9){\line(1,-1){6}}

\put(1.7,3.8){\makebox(1,0.5)[l]{\small $1$}}

\put(7.4,3.5){\makebox(1,0.5)[l]{\small $1$}}
\put(9.4,3.5){\makebox(1,0.5)[l]{\small $1$}}
\put(7.4,5.5){\makebox(1,0.5)[l]{\small $1$}}
\put(8.1,4.2){\makebox(1,0.5)[l]{\small $2$}}

\put(15.4,3.5){\makebox(1,0.5)[l]{\small $1$}}
\put(17.4,3.5){\makebox(1,0.5)[l]{\small $1$}}
\put(19.4,3.5){\makebox(1,0.5)[l]{\small $1$}}
\put(16.1,4.2){\makebox(1,0.5)[l]{\small $2$}}
\put(18.1,4.2){\makebox(1,0.5)[l]{\small $2$}}
\put(15.4,5.5){\makebox(1,0.5)[l]{\small $1$}}
\put(17.4,5.5){\makebox(1,0.5)[l]{\small $1$}}
\put(16.1,6.2){\makebox(1,0.5)[l]{\small $2$}}
\put(15.4,7.5){\makebox(1,0.5)[l]{\small $1$}}

\put(1,1){\makebox(1,0.5)[l]{\small (a)  $T_1$}}
\put(8,1){\makebox(1,0.5)[l]{\small (b) $T_2$}}
\put(17,1){\makebox(1,0.5)[l]{\small (c) $T_3$}}

\put(3,-1){\makebox(1,0.5)[l]{\small Figure 5. The levels of faces
for $T_1,T_2,T_3$.}}
\end{picture}
\end{center}

Assume that $n\geq 4$ and the assertion holds for smaller $n$. We
delete the boundary of the outer faces from $T_n$ (the vertices and
the edges on this boundary are deleted). Then the resulting graph is
$T_{n-3}$. In this transformation, all faces with levels $1$ and $2$
are removed. Therefore $\lambda_{\max}(T_n)=\lambda_{\max}(T_{n-3})
+2$. By induction hypothesis, we have
$$\lambda_{\max}(T_n)=\left\lceil\frac{2(n-3)}{3}\right\rceil+2=
\left\lceil\frac{2n}{3}\right\rceil.$$ For example,
$\lambda_{\max}(T_4)=\lambda_{\max}(T_1)+2=3$, as shown in Figure
4(a).

We next show the lower bound
\begin{equation}
\sigma_T(G)\geq \lambda_{\max}+1.
\end{equation}
In fact, let $T$ be any given spanning tree of $G$. Then the cotree
$\bar T$ determines a spanning tree $\bar T^*$ in $G^*$. Similar to
the previous case, suppose that $f_0$ is a face with the maximum
level $\lambda=\lambda_{\max}$. Then the distance between $f_0$ and
$O$ in $\bar T^*$ is at least $\lambda$. Let $P^*$ be the path from
$f_0$ to $O$ in $\bar T^*$ with the last edge $e^*_0$ incident with
$O$. The tree-edge $e^*_0$ in the spanning tree $\bar T^*$
determines a fundamental edge-cut $C^*$, which corresponds to a
fundamental cycle $C$ with respect to $T$ in $G$. This fundamental
cycle $C$ is determined by the cotree edge $e_0$ on the boundary of
the outer face that corresponds to the edge $e^*_0$ in $P^*$.
Without loss of generality, assume that $e_0$ is on the horizontal
line $R_0=\{(x,y)\in {\mathbf R}^2: y=0\}$. Let $P^0$ be the
shortest path from $f_0$ to $O$ passing though $R_0$ in $G^*$.
Suppose that $e'_0$ is the edge in $R_0$ which corresponds the last
edge of $P^0$. Denote by $\partial(P^0)$ the boundary of the region
composed of the $\lambda$ faces of $P^0$. Then $\partial(P^0)$ is a
$(1\times \tfrac 12\lambda)$ rectangle (if $\lambda$ is even) or a
$(1\times \tfrac 12(\lambda-1))$ rectangle plus a triangle of $f_0$
at the top (if $\lambda$ is odd). It can be seen that the triangle
face at the top (or at the bottom) of $\partial(P^0)$ has two
boundary edges, and each of the other triangle faces has one
boundary edge. Hence the length of $\partial(P^0)$ is $\lambda+2$.
We draw $\lceil \tfrac12 \lambda\rceil$ horizontal straight lines
passing through the midpoints of the boundary edges in
$\partial(P^0)$. Then each of these straight lines intersects the
cycle $C$ twice. When $\lambda$ is even, the $\tfrac 12\lambda$
straight lines intersect the cycle $C$ at $\lambda$ edges. Besides,
$C$ must have at least two more horizontal edges (one is $e_0$ and
one in $f_0$). Hence the length of $C$ is at least $\lambda+2$. When
$\lambda$ is odd, the $\tfrac 12(\lambda+1)$ straight lines
intersect the cycle $C$ at $\lambda+1$ edges. And $C$ has one more
horizontal edge $e_0$. Thus the length of $C$ is at least
$\lambda+2$. Therefore, for any spanning tree $T$, we find a
fundamental cycle $C$ with length at least $\lambda+2$. By the
arbitrariness of $T$, the above lower bound (4) is proved.

Conversely, we can construct an optimal spanning tree $T$ as follows:\\
{\indent} (1) Take a face $f_0$ with the maximum level
$\lambda_{\max}$. \\
{\indent} (2) Take the horizonal line $H$ containing the horizontal
edge of $f_0$, and take the vertical line $V$ containing the
vertical edge of $f_0$.\\
{\indent} (3) In the part below $H$, take every vertical line
intersecting $H$; In the part above $H$, take every horizontal line
intersecting $V$. \\
{\indent} (4) In the remaining part of the lower right corner, take
all horizontal lines; In the remaining part of the upper left
corner, take all vertical lines. An example of $T_4$ can be seen in
Figure 4(b). It is easy to check that this spanning tree attain the
above lower bound. This completes the proof. $\Box$

\subsection{Triangulated-rectangular grids}
\hspace*{0.5cm} Finally, we consider the triangulated-rectangular
grids by the same method. A triangulated-rectangular grid $T_{m,n}$
is defined as follows: the vertex set $V(T_{m,n})$ is $\{(x,y)\in
{\mathbf Z}^2: 0\leq y\leq m-1, 0\leq x\leq n-1\}$, and two vertices
$(x,y)$ and $(x',y')$ are joined by an edge if $|x-x'|+|y-y'|=1$ or
$|x-x'|+|y-y'|=2$ and $x+y=x'+y'$, as shown in Figure 6. Clearly,
$T_{m,n}$ can be obtained from the rectangular grids $P_m\times P_n$
by adding slant edges.

\begin{center}
\setlength{\unitlength}{0.45cm}
\begin{picture}(12,11)

\multiput(1,3)(2,0){6}{\circle*{0.3}}
\multiput(1,5)(2,0){6}{\circle*{0.3}}
\multiput(1,7)(2,0){6}{\circle*{0.3}}
\multiput(1,9)(2,0){6}{\circle*{0.3}}
\multiput(1,11)(2,0){6}{\circle*{0.3}}

\put(1,3){\line(1,0){10}} \put(1,5){\line(1,0){10}}
\put(1,7){\line(1,0){10}} \put(1,9){\line(1,0){10}}
\put(1,11){\line(1,0){10}} \put(1,3){\line(0,1){8}}
\put(3,3){\line(0,1){8}} \put(5,3){\line(0,1){8}}
\put(7,3){\line(0,1){8}} \put(9,3){\line(0,1){8}}
\put(11,3){\line(0,1){8}} \put(1,5){\line(1,-1){2}}
\put(1,7){\line(1,-1){4}} \put(1,9){\line(1,-1){6}}
\put(1,11){\line(1,-1){8}} \put(3,11){\line(1,-1){8}}
\put(5,11){\line(1,-1){6}} \put(7,11){\line(1,-1){4}}
\put(9,11){\line(1,-1){2}}

\put(1.4,3.4){\makebox(1,0.5)[l]{\small $1$}}
\put(3.4,3.4){\makebox(1,0.5)[l]{\small $1$}}
\put(5.4,3.4){\makebox(1,0.5)[l]{\small $1$}}
\put(7.4,3.4){\makebox(1,0.5)[l]{\small $1$}}
\put(9.4,3.4){\makebox(1,0.5)[l]{\small $1$}}
\put(1.4,5.4){\makebox(1,0.5)[l]{\small $1$}}
\put(1.4,7.4){\makebox(1,0.5)[l]{\small $1$}}
\put(1.4,9.4){\makebox(1,0.5)[l]{\small $1$}}
\put(2.1,10.1){\makebox(1,0.5)[l]{\small $1$}}
\put(4.1,10.1){\makebox(1,0.5)[l]{\small $1$}}
\put(6.1,10.1){\makebox(1,0.5)[l]{\small $1$}}
\put(8.1,10.1){\makebox(1,0.5)[l]{\small $1$}}
\put(10.1,10.1){\makebox(1,0.5)[l]{\small $1$}}
\put(10.1,8.1){\makebox(1,0.5)[l]{\small $1$}}
\put(10.1,6.1){\makebox(1,0.5)[l]{\small $1$}}
\put(10.1,4.1){\makebox(1,0.5)[l]{\small $1$}}

\put(2.1,4.1){\makebox(1,0.5)[l]{\small $2$}}
\put(4.1,4.1){\makebox(1,0.5)[l]{\small $2$}}
\put(6.1,4.1){\makebox(1,0.5)[l]{\small $2$}}
\put(8.1,4.1){\makebox(1,0.5)[l]{\small $2$}}
\put(2.1,6.1){\makebox(1,0.5)[l]{\small $2$}}
\put(2.1,8.1){\makebox(1,0.5)[l]{\small $2$}}
\put(3.4,9.4){\makebox(1,0.5)[l]{\small $2$}}
\put(5.4,9.4){\makebox(1,0.5)[l]{\small $2$}}
\put(7.4,9.4){\makebox(1,0.5)[l]{\small $2$}}
\put(9.4,9.4){\makebox(1,0.5)[l]{\small $2$}}
\put(9.4,7.4){\makebox(1,0.5)[l]{\small $2$}}
\put(9.4,5.4){\makebox(1,0.5)[l]{\small $2$}}

\put(3.4,5.4){\makebox(1,0.5)[l]{\small $3$}}
\put(5.4,5.4){\makebox(1,0.5)[l]{\small $3$}}
\put(7.4,5.4){\makebox(1,0.5)[l]{\small $3$}}
\put(3.4,7.4){\makebox(1,0.5)[l]{\small $3$}}
\put(4.1,8.1){\makebox(1,0.5)[l]{\small $3$}}
\put(6.1,8.1){\makebox(1,0.5)[l]{\small $3$}}
\put(8.1,8.1){\makebox(1,0.5)[l]{\small $3$}}
\put(8.1,6.1){\makebox(1,0.5)[l]{\small $3$}}

\put(4.1,6.1){\makebox(1,0.5)[l]{\small $4$}}
\put(6.1,6.1){\makebox(1,0.5)[l]{\small $4$}}
\put(5.4,7.4){\makebox(1,0.5)[l]{\small $4$}}
\put(7.4,7.4){\makebox(1,0.5)[l]{\small $4$}}

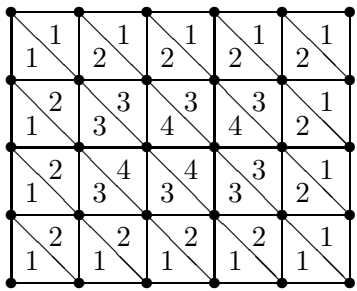
\put(-1,1){\makebox(1,0.5)[l]{\small Figure 6. Triangulated-
rectangle grid $T_{5,6}$}}
\end{picture}
\end{center}

{\bf Theorem 4.3.} \ For the triangulated-rectangular grids
$T_{m,n}$ with $2\leq m\leq n$, we have
$$\sigma_T(T_{m,n})=m.$$

{\bf Proof.} \ Let $G=T_{m,n}$ ($2\leq m\leq n$). We first claim
that
$$\lambda_{\max}(G)=m-1.$$
By induction on $m$. When $m=2$, all faces have level $1$, so
$\lambda_{\max}(G)=1$; When $m=3$, it is also evident that
$\lambda_{\max}(G)=2$. Assume that $m\geq 4$ and the claim holds for
smaller $m$. We delete the boundary of the outer faces from $G$, so
that the remaining graph is $G'=T_{m-2,n-2}$. In this
transformation, all faces with levels $1$ and $2$ are removed.
Therefore $\lambda_{\max}(G)=\lambda_{\max}(G')+2$. By induction
hypothesis, $\lambda_{\max}(G')=m-2-1=m-3$. Hence $\lambda_{\max}(G)
=m-3+2=m-1$ and the claim follows.

Moreover, by the same method of the previous case we obtain the
lower bound
\begin{equation}
\sigma_T(G)\geq \lambda_{\max}+1.
\end{equation}

Conversely, we can construct an optimal tree $T$ by taking all
columns and the edges from $(\lfloor m/2\rfloor,j)$ to $(\lceil m/2
\rceil, j+1)$ for $1\leq j\leq n-1$ (see Figure 7). It is easy to
check that this spanning tree attain the above lower bound. The
proof is complete. $\Box$

\begin{center}
\setlength{\unitlength}{0.45cm}
\begin{picture}(20,10)
\multiput(1,4)(2,0){4}{\circle*{0.3}}
\multiput(1,6)(2,0){4}{\circle*{0.3}}
\multiput(1,8)(2,0){4}{\circle*{0.3}}

\put(1,4){\line(0,1){4}} \put(3,4){\line(0,1){4}}
\put(5,4){\line(0,1){4}} \put(7,4){\line(0,1){4}}
\put(1,6){\line(1,0){6}}

\bezier{30}(1,4)(4,4)(7,4) \bezier{30}(1,8)(4,8)(7,8)
\bezier{15}(1,6)(2,5)(3,4) \bezier{30}(1,8)(3,6)(5,4)
\bezier{30}(3,8)(5,6)(7,4) \bezier{15}(5,8)(6,7)(7,6)

\multiput(12,4)(2,0){5}{\circle*{0.3}}
\multiput(12,6)(2,0){5}{\circle*{0.3}}
\multiput(12,8)(2,0){5}{\circle*{0.3}}
\multiput(12,10)(2,0){5}{\circle*{0.3}}

\put(12,4){\line(0,1){6}} \put(14,4){\line(0,1){6}}
\put(16,4){\line(0,1){6}} \put(18,4){\line(0,1){6}}
\put(20,4){\line(0,1){6}} \put(12,8){\line(1,-1){2}}
\put(14,8){\line(1,-1){2}} \put(16,8){\line(1,-1){2}}
\put(18,8){\line(1,-1){2}}

\bezier{30}(12,4)(16,4)(20,4) \bezier{30}(12,6)(16,6)(20,6)
\bezier{30}(12,8)(16,8)(20,8) \bezier{30}(12,10)(16,10)(20,10)
\bezier{10}(12,6)(13,5)(14,4) \bezier{10}(14,6)(15,5)(16,4)
\bezier{10}(12,10)(13,9)(14,8) \bezier{10}(16,6)(17,5)(18,4)
\bezier{10}(14,10)(15,9)(16,8) \bezier{10}(18,6)(19,5)(20,4)
\bezier{10}(16,10)(17,9)(18,8) \bezier{10}(18,10)(19,9)(20,8)

\put(1.5,2.2){\makebox(1,0.5)[l]{\small (a) $m$ is odd}}
\put(13,2.2){\makebox(1,0.5)[l]{\small (b) $m$ is even}}
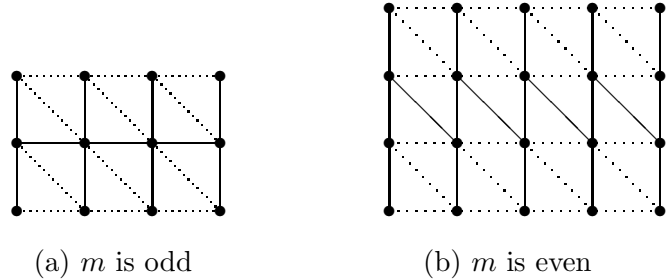
\put(4,0.5){\makebox(1,0.5)[l]{\small Figure 7. Optimal trees of
$T_{m,n}$.}}
\end{picture}
\end{center}


\begin{thebibliography}{abc}

\bibitem{Bondy08}  J.A. Bondy and U.S.R. Murty, Graph Theory.
Springer-Verlag, Berlin, 2008.

\bibitem{Bodlae11}  H.L. Bodlaender, K. Kozawa, T. Matsushima and
Y. Otachi, Spanning tree congestion of $k$-outerplanar graphs.
Discrete Mathematics, 311 (2011) 1040-1045.

\bibitem{Bodlae12} H.L. Bodlaender, F.V. Fomin, P.A. Golovach,
Y. Otachi and E.J. van Leeuwen, Parameterized complexity of the
spanning tree congestion problem. Algorithmica, 64 (2012) 85-111.

\bibitem{Brand04}  A. Brandst\"adt, F.F. Dragan, H.-O. Le  and
V.B. Le, Tree spanners on chordal graphs: complexity and algorithms.
Theoretical Computer Science, 310 (2004) 329-354.

\bibitem{Brand07} A. Brandst\"adt, F.F. Dragan, H.-O. Le,
V.B. Le and R. Uehara, Tree spanners for bipartite graphs and probe
interval graphs. Algorithmica, 47 (2007) 27-51.

\bibitem{Cai95} L. Cai and  D.G. Corneil, Tree spanners, SIAM
Journal of Discrete Mathematics, 8 (1995) 359-387.

\bibitem{Caste09}  A. Castej\'on and M.I. Ostrovskii, Minimum congestion
spanning trees of grids and discrete toruses, Discussiones
Mathematicae Graph Theory, 29 (2009) 511-519.

\bibitem{Chung88}  F.R.K. Chung, Labelings of graphs.
In: L.W. Beineke and R.J. Wilson (eds.), Selected Topics in Graph
Theory, Vol.3 (1988) 151-168.

\bibitem{Diaz02} J. Diaz, J. Petit and M. Serna, A survey of
graph layout problems. ACM Computing Surveys, 34 (2002) 313-356.

\bibitem{Fekete01} S.P. Fekete and J. Kremer, Tree spanners in
planar graphs. Discrete Applied Mathematics, 108 (2001) 85-103.

\bibitem{Fomin11}  F.V. Fomin, P.A. Golovach and E.J. van Leeuwen,
Spanners on bounded degree graphs. Information Processing Letters,
111 (2011) 142-144.

\bibitem{Galb03}  G. Galbiati, On finding cycle bases and fundamental
cycle bases with a shortest maximal cycles. Information Processing
Letters, 88 (2003) 155-159.

\bibitem{Golumbic80} M.C. Golumbic, Algorithmic Graph Theory and
Perfect graphs. Academic Press, New York 1980.

\bibitem{Hruska08}  S.W. Hruska, On tree congestion of graphs,
Discrete Mathematics, 308 (2008) 1801-1809.

\bibitem{Kozawa09} K. Kozawa, Y. Otachi and K. Yamazaki,  On
spanning tree congestion of graphs, Discrete Mathematics, 309 (2009)
4215-4224.

\bibitem{Lieb08} C. Liebchen and G. W\"unsch, The zoo of tree
spanner problems, Discrete Applied Mathematics, 156 (2008) 569-587.

\bibitem{Lin14} L. Lin,  Y. Lin and D. West, Cutwidth of triangular grids,
Discrete Mathematics, 331 (2014) 89-92.

\bibitem{Madan96} M.S. Madanlal,  G. Venkatesan and P.C. Rangan,
Tree 3-spanners on interval, permutation and regular bipartite
graphs. Information Processing Letters, 59 (1996) 97-102.

\bibitem{Okamoto11}  Y. Okamoto, Y. Otachi, R. Uehara and T. Uno,
Hardness results and an exact exponential algorithm for the spanning
tree congestion problem. Journal of Graph Algorithms and
Applications, 15(6) (2011) 727-751.

\bibitem{Ostrov04}  M.I. Ostrovskii, Minimal congestion trees.
Discrete Mathematics, 285 (2004) 219-326.

\bibitem{Ostrov10}  M.I. Ostrovskii, Minimum congestion spanning trees
in planar graphs. Discrete Mathematics, 310 (2010) 1204-1209.

\bibitem{Peleg89} D. Peleg and J.J. Ullman, An optimal
synchroniser for the hypercube. SIAM Journal of Computing, 18(4)
(1989) 740-747.

\bibitem{Venka97}  G. Venkatesan, U. Rotics, M.S. Madanlal, J.A. Makowsy
and C.P. Rangan,  Restrictions of minimum spanner problems.
Information and Computation, 136 (1997) 143-164.

\end{thebibliography}
\end{document}